\documentclass{article}
\usepackage{arxiv}
\usepackage[latin1,utf8]{inputenc} 
\usepackage{lmodern}
\usepackage{url}            
\usepackage{booktabs}       
\usepackage{amsmath}
\usepackage{amsfonts}		
\usepackage{amssymb}
\usepackage{mathtools}
\usepackage{commath}
\usepackage{nicefrac}       
\usepackage{microtype}      
\usepackage{graphicx}
\usepackage{cite}
\usepackage{doi}

\usepackage{hyperref} 
\hypersetup{
	colorlinks=true,
	linkcolor=blue,
	filecolor=magenta,
	urlcolor=cyan,
}

\usepackage{enumitem}
\usepackage{accents}
\usepackage{siunitx}
\usepackage{multirow}
\usepackage{array}
\newcolumntype{C}[1]{>{\centering\arraybackslash}m{#1}}
\usepackage{subfig}

\newcommand{\rh}{\tilde{r}}
\newcommand{\nh}{\tilde{n}}
\newcommand{\Nh}{\tilde{N}}

\newcommand{\Jhatxi}{J_{\hat{\xi}}}
\newcommand{\Jbrevexi}{J_{\breve{\xi}}}

\newcommand{\JbreveXi}{J_{\breve{\Xi}}}

\newcommand{\fcirc}{\underaccent{\circ}{\hat{f}}}

\newcommand{\res}{\text{res}}

\newcommand{\DIg}{\Dif_{\,1}\!g}
\newcommand{\DIIg}{\Dif_{\,2}\!g}

\newcommand{\DIres}{\Dif_{\,1}\!\text{res}(\hat{\xi},\hat{u};\varphi_\chi)\!\cdot\!\delta\hat{\xi}}
\newcommand{\DIIres}{\Dif_{\,2}\!\text{res}(\hat{\xi},\hat{u};\varphi_\chi)\!\cdot\!\delta\hat{u}}

\DeclareMathOperator{\trace}{tr}

\newcommand{\Ca}{\text{Ca}}

\title{Fully Implicit Spectral Boundary Integral Computation of Red Blood Cell Flow}

\author{
	\hspace{1mm}Pei Chuan Chao \\
	Department of Civil, Structural and Environmental Engineering \\
	University at Buffalo\\
	212 Ketter Hall, Buffalo, NY 14260 \\
	\texttt{peichuan@buffalo.edu} \\
	\And
	\href{https://orcid.org/0000-0003-3418-4883}
	{\includegraphics[scale=0.06]{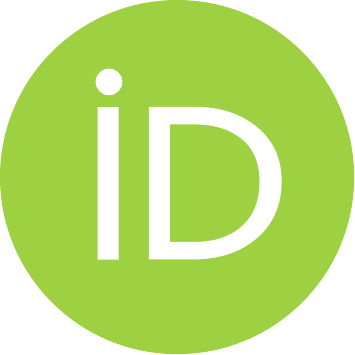}}
	 \hspace{1mm}Ali G\"urb\"uz \\
	 Department of Civil, Structural and Environmental Engineering \\
	 University at Buffalo\\
	 212 Ketter Hall, Buffalo, NY 14260 \\
	 \texttt{aligurbu@buffalo.edu} \\
	\And
	\href{https://orcid.org/0000-0002-4324-2240}
	{\includegraphics[scale=0.06]{orcid.pdf}}
	\hspace{1mm}Frederick Sachs \\
	Department of Physiology and Biophysics \\
	University at Buffalo\\
	304 Cary Hall, Buffalo, NY 14214 \\
	\texttt{sachs@buffalo.edu} \\
	\And
	\href{https://orcid.org/0000-0001-6661-5567}
	{\includegraphics[scale=0.06]{orcid.pdf}}
     \hspace{1mm}M.~V.~Sivaselvan 
     \thanks{Corresponding author} \\
	 Department of Civil, Structural and Environmental Engineering \\
	 University at Buffalo\\
	 212 Ketter Hall, Buffalo, NY 14260 \\
	 \texttt{mvs@buffalo.edu} \\
}

\hypersetup{
pdftitle={Fully Implicit Spectral Boundary Integral Computation of Red Blood Cell Flow},
pdfsubject={},
pdfauthor={Pei Chuan Chao, Ali G\"urb\"uz, Frederick Sachs, M.~V.~Sivaselvan},
pdfkeywords={Weakly singular boundary integral,
			 Fully implicit time integration,
			 Adaptive time stepping,
			 Spherical harmonics,
			 Elastic thin-shell mechanics,
			 Jacobian-vector multiplication},
}

\begin{document}
\maketitle

\begin{abstract}
	\noindent
	An approach is presented for implicit time integration in computations 
	of red blood cell flow by a spectral boundary integral method.
	The flow of a red cell in ambient fluid is represented as a boundary integral equation (BIE),
	whose structure is that of an implicit ordinary differential equation (IODE).
	The cell configuration and velocity field are discretized with spherical harmonics.
	The IODE is integrated in time using a multi-step implicit method based 
	on backward difference formulas, with variable order and adaptive time stepping
	controlled by local truncation error and convergence of Newton iterations.
	Jacobians of the IODE, required for Newton's method,
	are implemented as Jacobian matrix-vector products that are nothing but
	directional derivatives.
	Their computation is facilitated by the weakly singular format of the BIE,
	and these matrix-vector products themselves amount to computing a second BIE.
	Numerical examples show that larger time steps are possible,
	and that the number of matrix-vector products is comparable to explicit methods.
\end{abstract}

\section{Introduction}\label{sec:Intro}
Cell-level computational modeling of blood flow 
has been of great interest, because such simulations
offer insight into the rheological character of blood in microcapillaries,
aid in design of microfluidic devices for research and diagnosis,
and have potential to guide targeted drug delivery and treatment
\cite{Freund14,SecombReview2017,ComputationalCellMechBook}.
Our interest in modeling red blood cell (RBC) flow stems
from seeking to explore and test hypotheses about the role of the PIEZO1 
mechanosensitive ion channel (MSC) in normal functioning of RBC;
evidence of such a role is offered by clinical studies 
and by electrophysiological and microfluidic experiments 
\cite{Bae2011,Zarychanski2012,Bae2013,WanPIEZO1ATP2015,DanielczokPIEZO2017}.
This short note describes a byproduct of this exploration,
an implicit time integration method for simulation of RBC flow.
Investigating the role of PIEZOs involves simulation of RBC flow in
microcapillaries of diameter less than the RBC size \cite{Gurbuz2020}.
In this note, however, the simpler computational problem of deformation of a single RBC
in ambient flow in an unbounded domain is considered,
to present the formulation and implementation of the implicit method.

RBC flow at the cellular level is typically modeled using Stokes equation
to represent the intra- and extracellular plasma,
with the cell membrane described by elastic thin shell mechanics.
Within this modeling framework, different mathematical formulations
and corresponding numerical methods have been adopted 
(see \cite{Freund14} for a comprehensive review).
Among these, a commonly employed approach is to formulate
the Stokes flow of the intra- and extracellular plasma as 
boundary integral equations (BIE) \cite{Pozrikidis92,Pozrikidis03,Pozrikidis05}.
Since RBC flow involves moving boundaries corresponding to the 
cell membranes, the BIE formulation, where different fields are computed 
only on boundaries, and not in the interiors, is beneficial. 
Further within the BIE setting, different schemes 
have been employed to discretize fields on the cell membrane, 
such as finite elements \cite{DoddiBagchi2009,Boedec2011}
and Fourier-like orthogonal spherical harmonic basis functions \cite{Zhao10,Veerapaneni11}.
In this short note, we use a BIE formulation and spherical harmonic basis functions
to discretize the RBC deformed geometry and the velocity field on the cell membrane.
The cell membrane is nearly inextensible, and this is represented 
in elastic shell models either as a direct constraint \cite{Veerapaneni11,Rahimian15},
or by a penalty approach 
by means of a sufficiently large membrane dilatation modulus \cite{Zhao10,FreundSpleen13}.
In this note, the latter approach is adopted.

The BIE for the cell membrane velocity, together with the membrane shell mechanics, constitutes an implicit ordinary differential equation (IODE) \cite{AscherPetzoldBook} in the membrane configuration (see equation \eqref{eq:IODE})\footnote{When a vessel is present as well in the model, it is an index-1 differential algebraic equation (DAE) \cite{AscherPetzoldBook}, with the traction and velocity fields on the vessel wall being algebraic variables; if the vessel were also deformable, the BIE system would again be an IODE. This is the case when cell membrane inextensibility is modeled using a penalty approach; when inextensibility is imposed as a constraint as in \cite{Veerapaneni11,Rahimian15}, then the BIE is a DAE.}.
Most commonly, this IODE is integrated in time using explicit schemes
such as forward Euler \cite{Ramanujan98,Pozrikidis03,Pozrikidis05,Zhao10,Zhao11}.
Explicit schemes have been successfully used even in 
challenging scenarios, for example simulating red cells passing through
very small capillaries and slits, and taking on complex configurations
consisting of high buckled modes of deformation \cite{Zhao10,FreundSpleen13}.
However, it has been demonstrated in \cite{Veerapaneni11} that the BIE IODE/DAE system
is stiff due to inextensibility of the membrane and 
higher order shape derivatives that occur in computing bending deformation.
Hence when explicit schemes are used, time step is restricted by stability
rather than by accuracy.
Furthermore, modeling the kinetics of the PIEZO1 MSC will
introduce additional time scales \cite{Bae2011,Bae2013}.
Therefore, robustness with respect to time step size is desirable,
and we seek to implement an implicit integration scheme.

Implicit methods for time integration of BIE arising in 
interface evolution problems such as RBC flow 
have been explored less in the literature. 
Dimitrakopoulos \cite{Dimitrakopoulos07} presented an implicit method 
to compute the evolution of a three-dimensional droplet in ambient fluid of a different viscosity.
The interface was modeled as having constant surface tension
(unlike more general elastic shell mechanics in \cite{Zhao10,FreundSpleen13,Veerapaneni11,Rahimian15} and here).
Special characteristics of the droplet mechanics and a perturbation approach
were used to reduce the nonlinear algebraic equation in each time step
to one in a scalar amplitude field of a prescribed vector field search direction.
A single-step implicit Runge-Kutta method and a multi-step method based on backward difference
formulas (BDF) were both considered for time discretization.
Veerapaneni et al \cite{Veerapaneni11}, 
building on their work on axisymmetric flows \cite{veerapaneni09numerical},
proposed a semi-implicit scheme for simulation of three-dimensional vesicle flows.
The fluids inside and the outside the vesicles had the same viscosity,
and so the BIE contained only the single layer potential. 
Inextensibility of the membrane was modeled as a constraint,
and the membrane mechanics was derived from a strain energy function of mean curvature.
Backward Euler discretization was used, and 
the semi-implicit (or linearly implicit) scheme consisted of splitting the computation of the
bending component of the membrane force in such a way that the
new configuration can be computed by solving a linear equation.
Rahimian et al \cite{Rahimian15} extended this work to account for 
different viscosities inside and outside the vesicle, so that
the BIE also included a double layer potential. 
The BIE was treated using a Galerkin approach (as also done in this note). 
The implicit-explicit split was again applied to the bending force,
to obtain a linearly implicit scheme.
A hierarchy of iterative methods were considered for linear system solution,
and their performances compared. Despite reduced overall computational cost 
compared to the explicit scheme, large numbers of iterations per time step were needed.
Quaife and Biros \cite{Quaife16} presented an adaptive semi-implicit method 
for two-dimensional vesicle suspensions.
By adopting a spectral deferred correction method, 
the simulation was allowed to automatically select time step size 
and method accuracy up to third order.

In this note, an implicit method is presented that is quite simple to implement.
The flow is formulated as a BIE and the membrane mechanics using
elastic thin shell theory, 
with in-plane inextensibility of the membrane captured by a penalty approach.
Spherical harmonic basis functions are used to discretize different fields 
on the cell membrane, and as in \cite{Zhao10},
the package SPHEREPACK \cite{adams1999spherepack} is used to 
perform computations related to spherical harmonics.
The BIE, including the stress kernel,
is set up in a \emph{weakly singular} format following \cite{liu91}, 
facilitating its linearization.
The IODE structure of the BIE is deliberately recognized, and
an off-the-shelf IODE/DAE solver \cite{Hindmarsh05} is used.
The solver is passed two problem-specific routines
--- one that evaluates the BIE producing a residual,
and another that computes Jacobian matrix-vector products.
Computing a Jacobian matrix-vector product 
amounts to evaluating another BIE (see equation \eqref{eq:DIres}),
and the matrix itself is not constructed.
In forming the Jacobian matrix-vector product, 
both the BIE kernels and the membrane force are linearized.
The IODE solver uses BDF for time discretization,
and adapts time step as well as BDF order based on truncation error.
The solver also balances truncation error and accuracy of solution
of the nonlinear algebraic equations in each time step,
resulting in a comparable number of Jacobian matrix-vector products
per time step as needed in an iterative solution of the linear system
arising in explicit methods.

This short note is organized as follows. 
The BIE for RBC flow is presented in Section \ref{sec:BIE}
in weakly singular form.
Section \ref{sec:shell} summarizes key results from shell mechanics.
The IODE structure of the BIE is discussed in Section \ref{sec:IODE},
which concludes by stating the two functions
--- BIE evaluation and Jacobian matrix-vector product computation,
to be provided to the IODE solver.
The BIE directional derivatives needed for implicit solution
are presented in Section \ref{sec:DirectionalDeriv}.
Numerical examples are presented in Section \ref{sec:Examples},
highlighting various characteristics of the implicit solution.
Finally some key observations and potential for further development
of this implicit scheme implementation are summarized in Section \ref{sec:Conclusion}.

\section{Boundary integral equation (BIE)} \label{sec:BIE}
\begin{figure}
	\centering
	\resizebox{\textwidth}{!}{\includegraphics{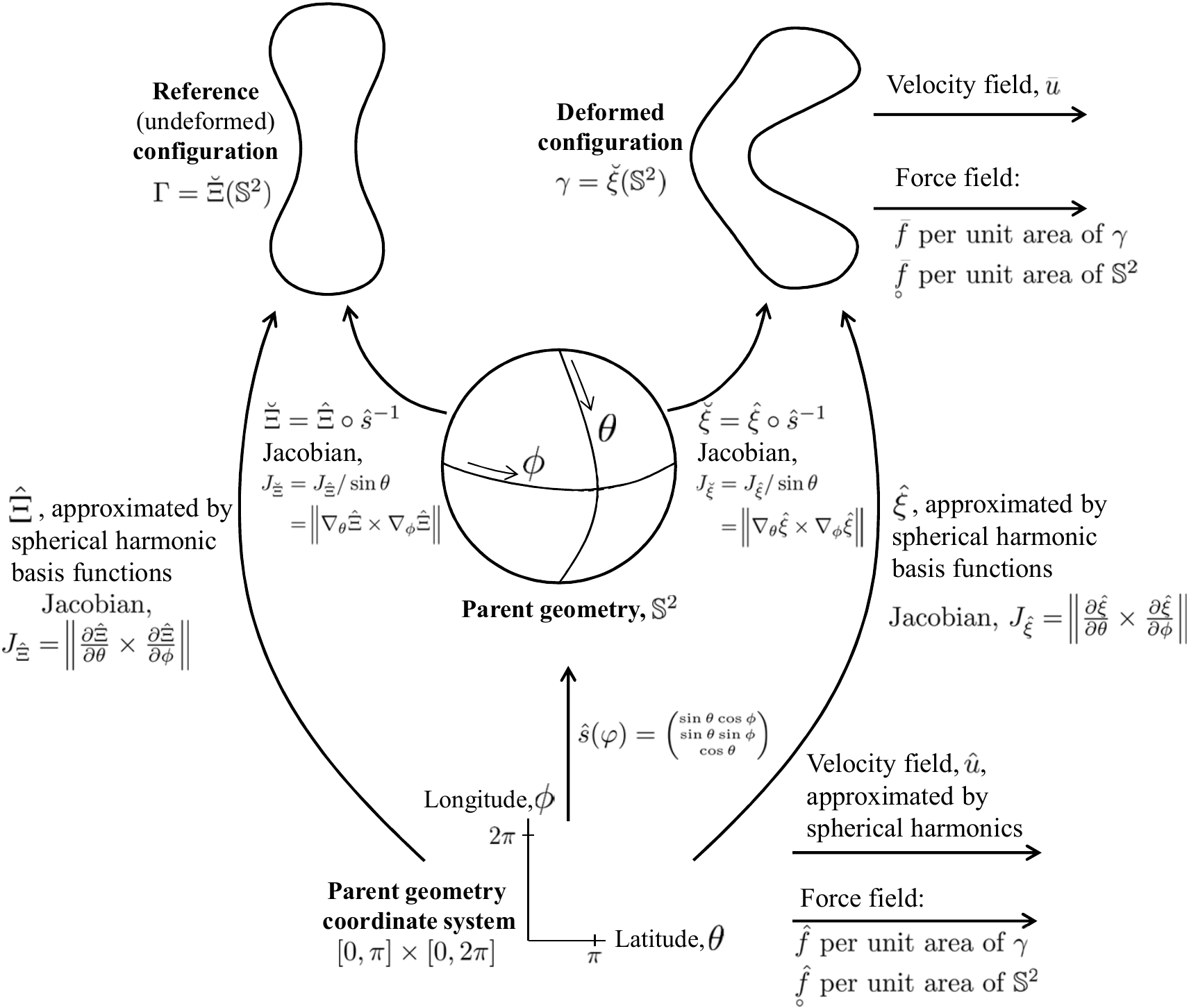}}
	\caption{Reference and deformed geometries of the RBC,
		the shapes of which are isomorphisms of the unit sphere (parent geometry);
		also shown is the coordinate system on the unit sphere.
		The hatted versions of the functions are the ones 
		approximated by spherical harmonic basis functions}
	\label{fig:Geometry}
\end{figure}
With the terminology and notation defined in Figure \ref{fig:Geometry},
the velocity field, $\bar{u}$, on the surface, $\gamma$, 
of a red cell in a three-dimensional ambient flow, $\bar{u}^\infty$, 
satisfies the following BIE
\begin{equation}
	\frac{\lambda+1}{2}\bar{u}(\chi) 
	+ \frac{\lambda-1}{8\pi}\!\int_\gamma K(\chi,\xi)\bar{u}(\xi)\dif\gamma \\
	= \bar{u}^\infty(\chi) 
	- \frac{1}{8\pi\mu}\!\int_\gamma G(\chi,\xi) \bar{f}(\xi)\dif\gamma 
	\;\; \forall \chi\!\in\!\gamma
	\label{eq:BIE}    
\end{equation}
where $\mu$ is the viscosity of the extracellular plasma,
$\lambda\mu$ that of the intracellular plasma,
$\bar{f}$ is the net force on the membrane applied by the fluids
on either side of it per unit area of the surface $\gamma$,
and the kernels,
\begin{equation}
	G(\chi,\xi) = \frac{1}{\norm{r}} (I + \rh \rh^\top); \quad 
	K(\chi,\xi) = \frac{6}{\norm{r}^2}(\rh^\top \nh)(\rh\rh^\top)
\end{equation}
are fundamental solutions for Stokes flow,
with $r = \chi - \xi$ and the unit vector, $\rh =\frac{r}{\norm{r}}$,
$\nh$ being the unit normal vector field on $\gamma$ pointing outward to ambient fluid.
Equation \eqref{eq:BIE} is the same as equation (11) in \cite{Pozrikidis03} 
(see also equation (2) on page 143 of \cite{Pozrikidis92}), 
equation (3a) in \cite{Rahimian15}, and 
equation (6) in \cite{Zhao10}, with appropriate mapping of notation.

Since the kernel $G$ contains $\norm{r}$ in the denominator,
the integral on the right hand side of equation \eqref{eq:BIE} is weakly singular,
whereas the integral on the left side is strongly singular
because the kernel $K$ contains $\norm{r}^2$ in the denominator \cite{Pozrikidis92}.
However, a weakly singular representation of the integral on the left may
be obtained by \emph{singularity subtraction} \cite{liu91};
subtracting and adding $\left[\int_\gamma K(\chi,\xi)\dif\gamma\right]\!\bar{u}(\chi)$,
and using the identity, 
$\int_\gamma K(\chi,\xi)\dif\gamma\!=\!-4\pi I$ \cite{Pozrikidis92,liu91}, we obtain
\begin{equation}
	\bar{u}(\chi) 
	+ \frac{\lambda-1}{8\pi}\!\int_\gamma K(\chi,\xi)(\bar{u}(\xi)-\bar{u}(\chi))\dif\gamma
	= \bar{u}^\infty(\chi) 
	- \frac{1}{8\pi\mu}\!\int_\gamma G(\chi,\xi) \bar{f}(\xi)\dif\gamma 
	\;\; \forall \chi\!\in\!\gamma
	\label{eq:BIE-weak}    
\end{equation}
where now, the integral on the left is also weakly singular,
and $(\bar{u}(\chi)-\bar{u}(\xi))$ may be thought of as ``cancelling'' 
one of the $\norm{r}$s in the denominator of $K$ as $\chi\xrightarrow{}\xi$.
This form allows using the method of \cite{GaneshGraham2004},
involving rotation of the spherical coordinate system to align the north pole with the source point, $\chi$,
the eigenfunction property and the additional theorem of spherical harmonics,
to compute both integrals
(see \cite{Gurbuz2020} for further details).
Furthermore, when the BIE is linearized in equation \eqref{eq:DIres} below,
the resulting integrals remain weakly singular,
and the method of \cite{GaneshGraham2004} can be used to compute them as well.

Introducing spherical coordinates, we have
\begin{equation}
	\begin{gathered}
		\hat{u}(\varphi_\chi) 
		+ \frac{\lambda-1}{8\pi}\!
		\int_{\phi=0}^{2\pi}\!\int_{\theta=0}^{\pi}
		K\left(\hat{\xi}(\varphi_\chi),\hat{\xi}(\varphi)\right)
		(\hat{u}(\varphi)-\hat{u}(\varphi_\chi))
		\Jhatxi(\varphi)
		\dif\theta\dif\phi \\
		= \hat{u}^\infty(\varphi_\chi) 
		- \frac{1}{8\pi\mu}\!
		\int_{\phi=0}^{2\pi}\!\int_{\theta=0}^{\pi} 
		G\left(\hat{\xi}(\varphi_\chi),\hat{\xi}(\varphi)\right) 
		\fcirc(\varphi) \sin\theta
		\dif\theta\dif\phi \\
		\;\, \forall \varphi_\chi\!\in\![0,\pi]\times[0,2\pi]
	\end{gathered}
	\label{eq:BIE-coord}    
\end{equation}
where $\varphi_\chi$ denotes the spherical coordinates 
that map to the source point, $\chi$.
Note that the force on the membrane in equation \eqref{eq:BIE-coord} is
expressed per unit area of the unit sphere parent geometry, $\fcirc$,
since this can be computed readily by equation \eqref{eq:Fcirc} below.
Moreover, since the Jacobian $J_{\hat{\xi}}$ does not appear in the 
right hand side integral of equation \eqref{eq:BIE-coord},
the linearization in equation \eqref{eq:DIres} below has a simpler form.

\section{Summary of elastic thin-shell mechanics results}\label{sec:shell}
This section summarizes the computation of the membrane force.
The formulas are completely equivalent to standard ones in the literature 
\cite{Flugge1972,Niordson1985,Zhao10}. However, they are arranged 
so that the force can be computed using 
only the operations of the gradient of scalar functions
and the divergence of vector fields on the unit sphere,
without invoking Christoffel symbols and derivatives of higher-order tensor fields.
The gradient and divergence are computed in spherical coordinates as
\begin{equation}
	\nabla\square 
	= \left( 
	\frac{\partial\square}{\partial\theta},
	\frac{1}{\sin\theta}\frac{\partial\square}{\partial\phi}
	\right)^\top; \quad
	\nabla\cdot\square 
	= \frac{1}{\sin\theta}
	\left( 
	\frac{\partial}{\partial\theta}(\square_\theta\sin\theta) +
	\frac{\partial}{\partial\phi}\square_\phi
	\right)
	\label{eq:GradDiv}
\end{equation}
directly using the \texttt{gradgc} and \texttt{divgc} functions in SPHEREPACK.
Note that the Jacobians in Figure \ref{fig:Geometry} and normal vector fields 
to the reference and deformed surfaces
can also be computed using the gradient operator,
\begin{equation}
	\begin{gathered}
		J_{\breve{\Xi}} = \norm{\nabla_\theta\hat{\Xi}\times\nabla_\phi\hat{\Xi}}; \quad
		\Nh(\varphi) = \frac{1}{J_{\breve{\Xi}}}\left(\nabla_\theta\hat{\Xi}\times\nabla_\phi\hat{\Xi}\right) \\
		J_{\breve{\xi}} = \norm{\nabla_\theta\hat{\xi}\times\nabla_\phi\hat{\xi}}; \quad
		\nh(\varphi) = \frac{1}{J_{\breve{\xi}}}\left(\nabla_\theta\hat{\xi}\times\nabla_\phi\hat{\xi}\right)
	\end{gathered}
	\label{eq:jacobian-normal}
\end{equation}

The force per unit area of $\mathbb{S}^2$ is computed by the following steps
(details can be found in \cite{Chao2020}).
\begin{enumerate}
	\item \textbf{Kinematics}: 
	The strain, $\epsilon$, and curvature change, $\kappa$, are computed as
	\begin{equation}
		\epsilon = \frac{1}{2}A^{-1}(a - A); \quad 
		\kappa = A^{-1}(b - B)
		\label{eq:StrainCurvature}
	\end{equation}
	where $A = \nabla\hat{\Xi}^\top\nabla\hat{\Xi}$ and $a = \nabla\hat{\xi}^\top\nabla\hat{\xi}$
	are matrices closely related to the metric tensor in the reference 
	and deformed configurations; 
	similarly, $B = \nabla\hat{\Xi}^\top\nabla\Nh$ and $b = \nabla\hat{\xi}^\top\nabla\nh$
	are related to the respective curvatures.
	
	\item \textbf{Constitutive equations}:
	The membrane in-plane force and moment matrices, $\mathfrak{n}$ and $\mathfrak{m}$,
	are then obtained using constitutive equations,
	\begin{equation}
		\begin{aligned}
			\mathfrak{n} &= E_\text{S}I 
			+ \left(\frac{E_\text{D}}{2}\log(\det(C))-E_\text{S}\right) C^{-\top} \\
			\mathfrak{m} &= E_\text{B} \kappa^\top
		\end{aligned}
		\label{eq:Constitutive}
	\end{equation}
	where $C = 2\epsilon + I$, and following the notation of \cite{Pozrikidis05,Zhao10},
	$E_\text{D}$ and $E_\text{S}$ are in-plane dilatational and shear moduli,
	and $E_\text{B}$ is a bending modulus. A neo-Hookean model is used for the 
	in-plane forces, and a linear model for the moment.
	
	\item \textbf{Equilibrium}: 
	The membrane out-of-plane shear, $\nh\mathfrak{q}$ is computed as
	\begin{equation}
		\nh\mathfrak{q} = \frac{1}{\JbreveXi\Jbrevexi}
		\left[\begin{array}{c|c}
			\nabla_\phi\hat{\xi}\times\mathcal{P}\mathfrak{r} &
			\mathcal{P}\mathfrak{r}\times\nabla_\theta\hat{\xi}
		\end{array}\right]
		\label{eq:out-of-plane-shear}
	\end{equation}
	where the normal projection, 
	$\mathcal{P} = I - \nh\nh^\top$,
	and $\mathfrak{r}=-\nabla\cdot\left(\nabla\hat{\xi}A^{-1}\mathfrak{m}J_{\breve{\Xi}}\right)$.
	Finally, the membrane force is calculated using
	\begin{equation}
		\fcirc(\varphi)
		= -\nabla \cdot 
		\left[
		\left(
		\nabla\hat{\xi}A^{-1}\mathfrak{n}
		+ \nabla\nh A^{-1}\mathfrak{m}
		+ \nh\mathfrak{q}
		\right)J_{\breve{\Xi}}
		\right]
		\label{eq:Fcirc}
	\end{equation}
\end{enumerate}

\section{Galerkin approach and implicit ordinary differential equation (IODE) structure}\label{sec:IODE}
\begin{figure}
	\centering
	\resizebox{0.95\textwidth}{!}{\includegraphics{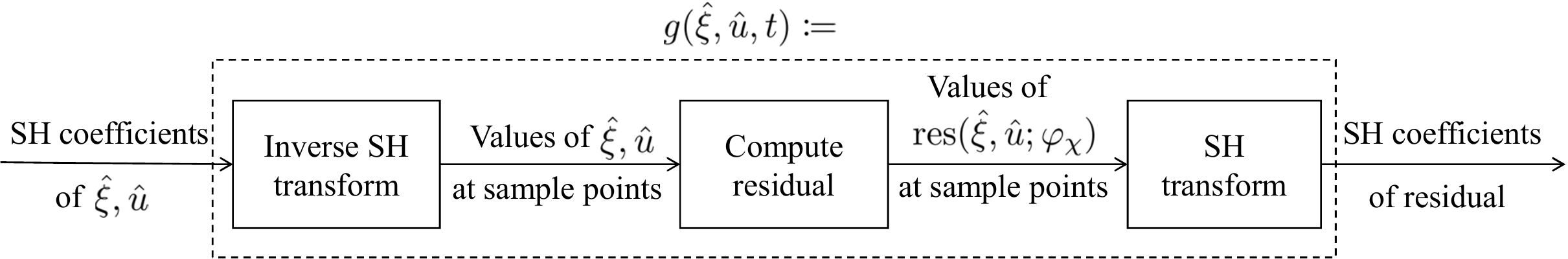}}
	\caption{Implicit ODE structure of BIE}
	\label{fig:IDOE-g}
\end{figure}

The residual of equation \eqref{eq:BIE-coord}, which we call $\res(\hat{\xi},\hat{u};\chi)$,
may be computed for any source point, $\chi$.
The BIE is satisfied if this residual is zero for all source points.
Instead of enforcing the BIE at selected source points,
the residual is computed using all the sample points as source points, 
and the spherical harmonic transform of the residual is set to zero.
This amounts to a Galerkin approach with the test functions being
the spherical harmonic basis functions (see \cite{Gurbuz2020} for details).
The process is depicted in Figure \ref{fig:IDOE-g}, 
and may be thought of as the equation,
\begin{equation}
	g(\hat{\xi}, \hat{u}, t) = 0
	\label{eq:IODE}
\end{equation}

Since the velocity field, $\hat{u}\!=\!\dot{\hat{\xi}}$, is
the rate of change of configuration, 
this is an implicit ordinary differential equation (IODE).

Numerical solutions to such IODE could be computed by means of 
a number of solvers such as 
\texttt{ode15i} in the MATLAB suite \cite{MATLAB:2019a},
\texttt{DASSL} \cite{dassl},
\texttt{DASPK} \cite{daspk} and
\texttt{IDAS} in the SUNDIALS suite \cite{Hindmarsh05}
that use multi-step methods based on backward difference formulas (BDF),
and \texttt{RADAU5} \cite{radau} based on an implicit Runge-Kutta method.
All such solvers require the derivatives, or Jacobians, $\DIg$ and $\DIIg$,
to use in different variations of Newton's method at each time-step.
Computing these Jacobian matrices is cumbersome and
computationally expensive.
However, as outlined below in Section \ref{sec:DirectionalDeriv},
the directional derivatives (or matrix-vector products),
$\DIg(\hat{\xi},\hat{u},t)\!\cdot\!\delta\hat{\xi}$ and
$\DIIg(\hat{\xi},\hat{u},t)\!\cdot\!\delta\hat{u}$
along $\delta\hat{\xi}$ and $\delta\hat{u}$ can be readily computed.
The solvers \texttt{DASPK} and \texttt{IDAS} allow 
a function that returns these matrix-vector products.
We use \texttt{IDAS} with its MATLAB interface.
This solver requires two functions representing the IODE,
\begin{enumerate}
	\item A function that given the state $(\hat{\xi},\hat{u})$ and time,
	computes the residual $g(\hat{\xi},\hat{u},t)$; 
	this is implemented as shown in Figure \ref{fig:IDOE-g}.
	\item A function that given the state, $(\hat{\xi},\hat{u})$, time,
	directions $(\delta\hat{\xi},\delta\hat{u})$, and a parameter, $\Lambda$,
	computes the sum of the directional derivatives,
	$\DIg(\hat{\xi},\hat{u},t)\!\cdot\!\delta\hat{\xi} +
	\Lambda\DIIg(\hat{\xi},\hat{u},t)\!\cdot\!\delta\hat{u}$.
	This is discussed in the next section.
\end{enumerate}

\section{Directional derivatives}\label{sec:DirectionalDeriv}
\begin{figure}
	\centering
	\resizebox{0.95\textwidth}{!}{\includegraphics{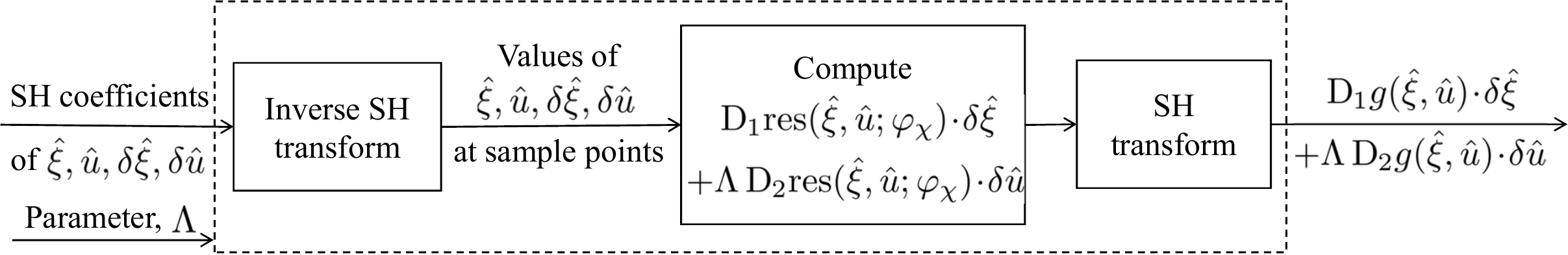}}
	\caption{Structure of function to compute directional derivatives 
		(Jacobian matrix-vector products) of the BIE IODE}
	\label{fig:IDOE-Dg}
\end{figure}

This section is on the computation of the directional derivatives 
of the IODE residual needed for implicit numerical solution.
The directional derivatives, $\DIres$ and $\DIIres$ residual of equation \eqref{eq:IODE}
are calculated, and required derivative of the IODE is obtained after a transformation
as shown in Figure \ref{fig:IDOE-Dg}.

The second of these directional derivatives is simple,
since the BIE \eqref{eq:BIE-coord} is linear in the velocity field,
and can be computed in the same manner as the BIE itself.
\begin{equation}
	\begin{gathered}
		\DIIres 
		= \delta\hat{u}(\varphi_\chi) \\
		+\frac{\lambda-1}{8\pi}\!
		\int_{\phi=0}^{2\pi}\!\int_{\theta=0}^{\pi}
		K\left(\hat{\xi}(\varphi_\chi),\hat{\xi}(\varphi)\right)
		(\delta\hat{u}(\varphi)-\delta\hat{u}(\varphi_\chi))
		\Jhatxi(\varphi)
		\dif\theta\dif\phi
	\end{gathered}
	\label{eq:DIIres}
\end{equation}
The derivative in a direction of configuration change 
requires more careful consideration.
\begin{equation}
	\begin{aligned}
		\DIres &= 
		\frac{\lambda-1}{8\pi}\!
		\int_{\phi=0}^{2\pi}\!\int_{\theta=0}^{\pi}
		\left[
		\delta\mathcal{K} \Delta\hat{u} \Jbrevexi
		+ \mathcal{K} \Delta\hat{u} \delta\Jbrevexi
		\right]
		\sin\theta\dif\theta\dif\phi \\
		&- \delta\hat{u}^\infty(\varphi_\chi)
		+ \frac{1}{8\pi\mu}\!
		\int_{\phi=0}^{2\pi}\!\int_{\theta=0}^{\pi} 
		\left[
		\delta\mathcal{G} \fcirc 
		+ \mathcal{G} \delta\fcirc 
		\right]
		\sin\theta\dif\theta\dif\phi
	\end{aligned}
	\label{eq:DIres}
\end{equation}
where we have used the abbreviations,
$\Delta\hat{u}\!=\!\hat{u}(\varphi)-\hat{u}(\varphi_\chi)$,
$\mathcal{K}\!=\!K(\hat{\xi}(\varphi_\chi),\hat{\xi}(\varphi))$ and
$\mathcal{G}\!=\!G(\hat{\xi}(\varphi_\chi),\hat{\xi}(\varphi))$,
and $\delta\square$ to represent the directional derivative, 
$\Dif\square(\hat{\xi})\!\cdot\!\delta\hat{\xi}$.
In the following, we define $\nu = \nabla_\theta\hat{\xi}\times\nabla_\phi\hat{\xi}$.
The directional derivatives in equation \eqref{eq:DIres},
obtained by straightforward but tedious calculation, are given by
\begin{equation}
	\begin{gathered}
		\delta r  = \delta\hat{\xi}(\varphi_\chi) - \delta\hat{\xi}(\varphi); \quad
		\delta\rh = \frac{1}{\norm{r}}(I-\rh\rh^\top)\delta r \\
		\delta\nu = \nabla_\theta\delta\hat{\xi}(\varphi)\times\nabla_\phi\hat{\xi}(\varphi) 
		+ \nabla_\theta\hat{\xi}(\varphi)\times\nabla_\phi\delta\hat{\xi}(\varphi) \\
		\delta\Jbrevexi = \nh^\top\delta\nu; \;\;
		\delta\nh = \frac{1}{\Jbrevexi}\mathcal{P}\delta\nu; \;\;
		\delta\mathcal{G}\fcirc
		= \frac{1}{\norm{r}}\left[ -(\rh^\top\delta r)\mathcal{G}\fcirc 
		+ (\rh^\top\fcirc)\delta\rh 
		+ (\delta\rh^\top\fcirc)\rh 
		\right] \\
		\begin{aligned}
			\delta\mathcal{K}\Delta\hat{u}
			= \frac{6}{\norm{r}^2}
			&\left[
			(\rh^\top\nh)(\rh^\top\Delta\hat{u})\left(\delta\rh-\frac{2}{\norm{r}}(\rh^\top\delta r)\rh\right)
			\right.\\
			&\left.
			+ \left( 
			(\rh^\top\Delta\hat{u})(\delta\rh^\top\nh + \rh^\top\delta\nh) 
			+ (\rh^\top\nh)(\delta\rh^\top\Delta\hat{u}) 
			\right)\rh  
			\vphantom{\left(\frac{2}{\norm{r}}\right)}  
			\right]
		\end{aligned}
	\end{gathered}
	\label{eq:DIres-DirectionalDeriv}
\end{equation}
The directional derivative of the membrane force is obtained 
by linearizing the shell mechanics relations from Section \ref{sec:shell}.
\begin{equation}
	\begin{aligned}
		\delta\epsilon &= \frac{1}{2}A^{-1}
		\left( 
		\nabla\delta\hat{\xi}^\top\nabla\hat{\xi} 
		+ \nabla\hat{\xi}^\top\nabla\delta\hat{\xi} 
		\right) \\
		\delta C &= 2\delta \epsilon \\
		\delta\kappa &= A^{-1}
		\left( 
		\nabla\delta\hat{\xi}^\top\nabla\nh 
		+ \nabla\hat{\xi}^\top\nabla\delta\nh
		\right) \\
		\delta\mathfrak{n} 
		&= \left[
		\left (
		\frac{E_\text{D}}{2}\trace(C^{-1}\delta C)
		\right)
		- \left(
		\frac{E_\text{D}}{2}\log(\det(C))-E_\text{S}
		\right)\delta C
		\right]C^{-\top} \\
		\delta\mathfrak{m} &= E_\text{B} \delta\kappa^\top \\
		\delta\mathfrak{r}
		&=-\nabla\cdot
		\left[
		\left(
		\nabla\delta\hat{\xi}A^{-1}\mathfrak{m}
		+ \nabla\hat{\xi}A^{-1}\delta\mathfrak{m}
		\right)\JbreveXi
		\right] \\
		\delta(\mathcal{P}\mathfrak{r})
		&= \mathcal{P}\delta\mathfrak{r}
		- (\nh^\top\mathfrak{r})\delta\nh
		- (\delta\nh^\top\mathfrak{r})\nh \\
		\delta(\nh\mathfrak{q})
		&= -\frac{\delta\Jbrevexi}{\Jbrevexi}(\nh\mathfrak{q})
		+\frac{1}{\JbreveXi\Jbrevexi}
		\left[\begin{array}{c|c}
			\!\nabla_\phi\delta\hat{\xi}\!\times\!\mathcal{P}\mathfrak{r}
			\!+\!\nabla_\phi\hat{\xi}\!\times\!\delta(\mathcal{P}\mathfrak{r}) &
			\!\mathcal{P}\mathfrak{r}\!\times\!\nabla_\theta\delta\hat{\xi}
			\!+\!\delta(\mathcal{P}\mathfrak{r})\!\times\!\nabla_\theta\hat{\xi}
		\end{array}\right] \\
		\delta \fcirc(\varphi)
		&= -\nabla\!\cdot\!
		\left[
		\left(
		\nabla\delta\hat{\xi}A^{-1}\mathfrak{n}
		+ \nabla\hat{\xi}A^{-1}\delta\mathfrak{n}
		+ \nabla\delta\nh A^{-1}\mathfrak{m}
		+ \nabla\nh A^{-1}\delta\mathfrak{m}
		+ \delta(\nh\mathfrak{q})
		\right)J_{\breve{\Xi}}
		\right]
	\end{aligned}
	\label{eq:Shell-linearization}
\end{equation}
Equation \eqref{eq:DIres} is itself a BIE with only weakly singular terms, 
and the integrals can again be computed using the method of \cite{GaneshGraham2004}.
The weakly singular representation of BIE \eqref{eq:BIE-coord}
facilitates the linearization in this convenient form.

\section{Numerical examples}\label{sec:Examples}
Four numerical examples are presented to demonstrate 
the performance of the proposed fully implicit method. 
Various aspects of performance such as 
(i) time steps used to achieve the specified tolerances, 
(ii) order of backward difference approximation used, 
(iii) number of function and derivative calls needed, 
(iv) conservation of different quantities, and
(v) comparison with the explicit time scheme, 
as well as the dependence of these on spatial discretization are illustrated. 
In each example, the motion of a red blood cell suspended in an ambient flow is considered. 
The four examples are 
\begin{enumerate}[label=(\alph*)]
	\item Shear flow ($k\!=\!100\SI{}{\per\second}$): An ambient shear flow parallel to the $x$ direction with a shear strain rate, $k =100\SI{}{\per\second}$ (as in \cite{Pozrikidis03, Veerapaneni11}), so that ambient velocity field, $u = (kz,0,0)^\top$, with an undeformed red cell placed initially making an angle of $45^{\circ}$ with the $x$ axis. Capillary number, $\Ca =\frac{\mu k a}{E_s} = 0.08$.
	\item Shear flow ($k\!=\!500\SI{}{\per\second}$): Same as (a), but with a greater shear rate, $k\!=\! 500\SI{}{\per\second}$ (as in \cite{Pozrikidis03, Veerapaneni11}), $\Ca = 0.4$.
	\item Parabolic flow: A parabolic ambient flow along the $x$ direction, $u_x = A(B - (y^2 + z^2))$, where A = 1 (\SI{35.46}{\per\micro\metre\per\second}) and B = 2.3 (\SI{18.3}{\micro\metre\squared}) are parameters that control the curvature and peak velocity of the flow, $\Ca = \frac{\mu A a^4}{E_B} = 14.95$.
	\item Relaxation: A highly deformed cell is suspended in ambient fluid at rest until it relaxes to its undeformed shape. A deformed shape from (b) is taken as the initial configuration.
\end{enumerate}

Following common practice in the literature \cite{Pozrikidis03, Pozrikidis05, Zhao10}, 
non-dimensional parameters are used. 
The equivalent cell radius, $a \approx\SI{2.82}{\micro\metre}$ is used as characteristic length, 
the inverse of a shear rate $k = 100$ $\SI{}{\per\second}$ as characteristic time, 
and plasma viscosity times this shear rate $\mu k  = \SI{0.12}{\Pa}$ as characteristic stress, 
where $\mu = \SI{1.2e-3}{\Pa\cdot\second}$. With this nondimensionalization, the cell membrane shear modulus, $E_\text{S} = 12.4$ $(4.2 \times 10^{-6}\text{N}/\text{m})$ and bending modulus, 
$E_\text{B} = 0.0669$ $(1.8 \times 10^{-19}\text{N}\cdot\text{m}$). The membrane dilatational modulus is set as a penalty parameter $E_\text{D} = 200$ $(6.8 \times 10^{-5}\text{N}/\text{m})$ to control membrane area change while preserving numerical conditioning. Viscosity ratio of the intracellular plasma to extracellular plasma, $\lambda = 5$. The ability of the membrane to deform under viscous stress is gauged by capillary number, $\Ca = \frac{\mu k a}{E_s}$, based on shear modulus for case (a) and (b), and $\Ca = \frac{\mu A a^4}{E_B}$, based on bending modulus for case (c) \cite{FreundSpleen13}. The initial undeformed shape of the cell is a biconcave disk, which following \cite{Pozrikidis03, Pozrikidis05} is described by
\begin{align}
	&\hat{\Xi}(\varphi) = 
	\left(
	\alpha \sin{\theta}\cos{\phi};
	\alpha \sin{\theta}\sin{\phi};
	\frac{\alpha}{2}(0.207+2.003\sin^2\theta-1.123\sin^4\theta)\cos{\theta}
	\right)^\top
\end{align}
where $\alpha\!=\!1.386$ is the non-dimensional maximum radius of the biconcave disk.

The relative tolerance for truncation error is set to $10^{-4}$ and absolute tolerance to $10^{-6}$ 
for the \texttt{IDAS} solver for all cases. A maximum time step restriction is set for each example
to minimize failed steps.

\begin{figure}
	\centering
	\begin{tabular}{C{0.6cm}C{2.1cm}C{11cm}}
		\hline
		&Flow type & Deformation snapshots\\
		\hline\\[-1em]
		(a)& Shear flow k = 100 \SI{}{\per\second} & \includegraphics[width=0.71\textwidth]{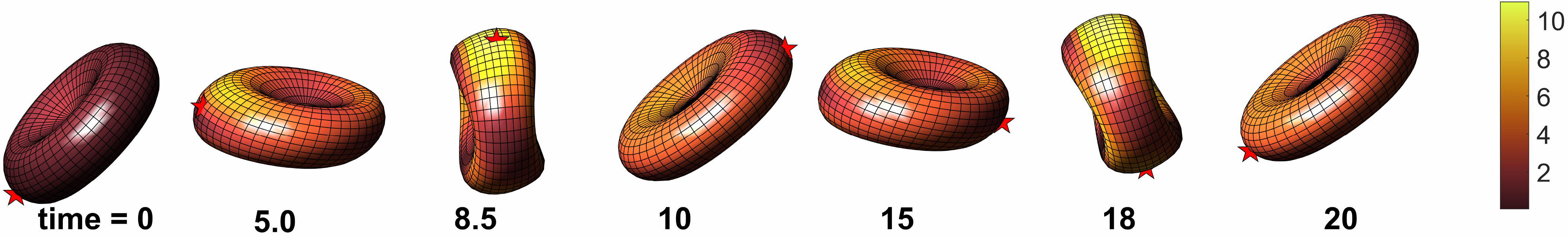} \\ 
		\hline\\[-1em]
		(b)& Shear flow k = 500 \SI{}{\per\second}  & \includegraphics[width=0.71\textwidth]{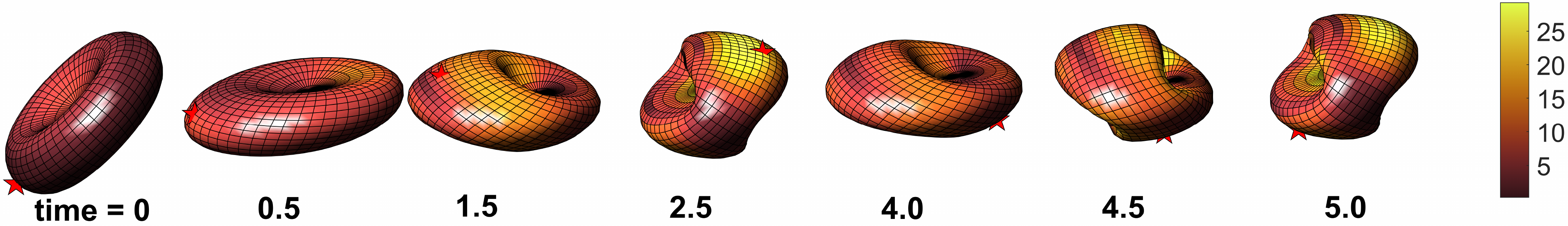} \\
		\hline\\[-1em]
		(c)& Parabolic flow & \includegraphics[width=0.71\textwidth]{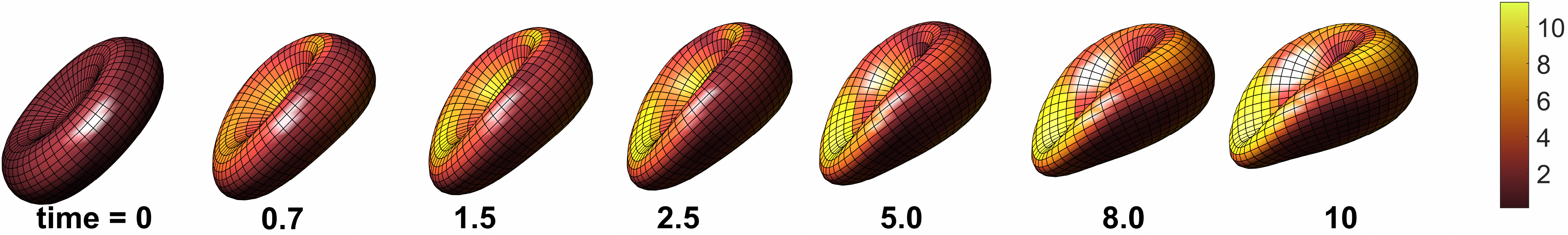} \\
		\hline\\[-1em]
		(d)& Relaxation & \includegraphics[width=0.71\textwidth]{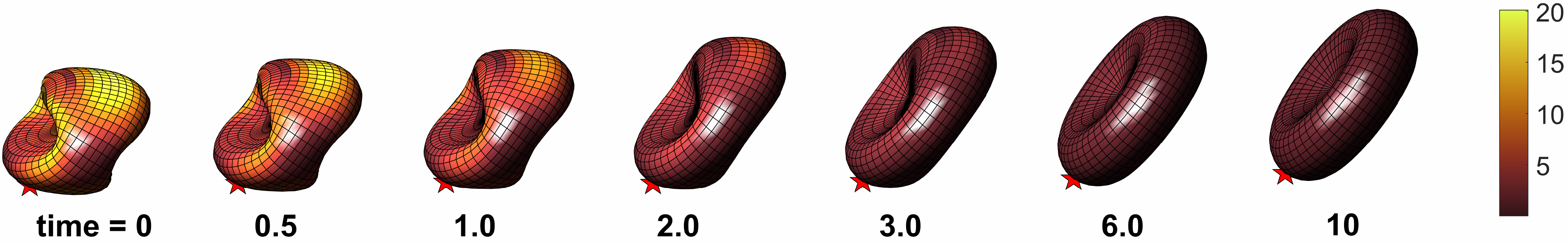}\\
		\hline
	\end{tabular}
	\caption{Evolution of cell deformation ($N$=24). Colors indicate the norm of the principal membrane stress.}
	\label{fig:deform}
\end{figure}

\begin{figure}
	\centering
	\begin{tabular}{C{0.6cm}C{2.1cm}C{11cm}}
		\hline
		&Flow type & Profile snapshots\\
		\hline\\[-1em]
		(a)& Shear flow k = 100 \SI{}{\per\second} & \includegraphics[width=0.73\textwidth]{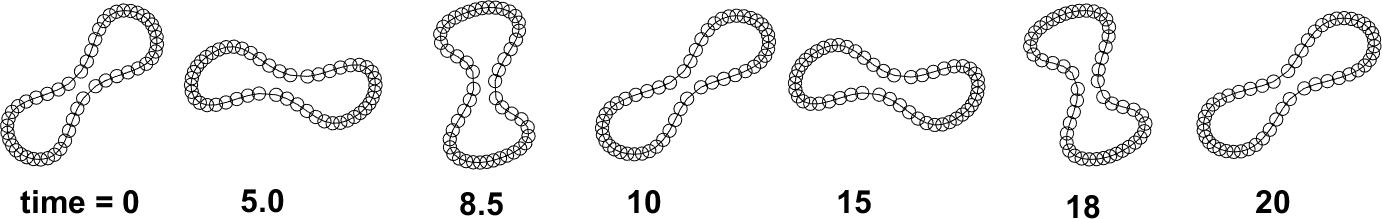} \\
		\hline\\[-1em]
		(b)& Shear flow k = 500 \SI{}{\per\second}  & \includegraphics[width=0.73\textwidth]{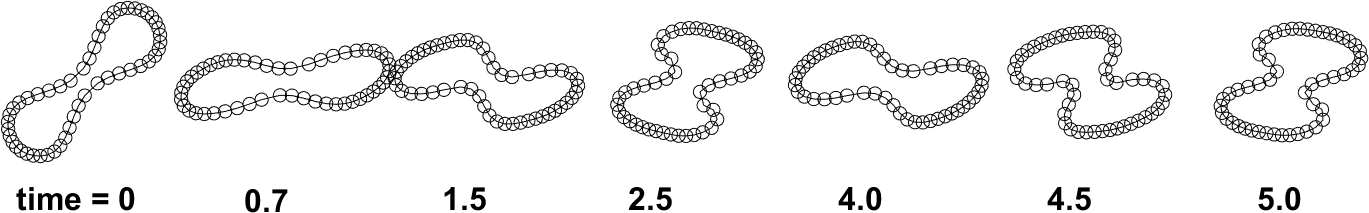} \\
		\hline\\[-1em]
		(c)& Parabolic flow & \includegraphics[width=0.73\textwidth]{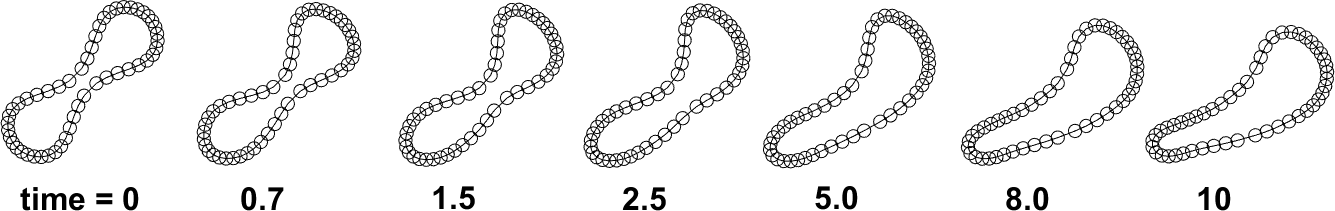} \\
		\hline\\[-1em]
		(d)& Relaxation & \includegraphics[width=0.73\textwidth]{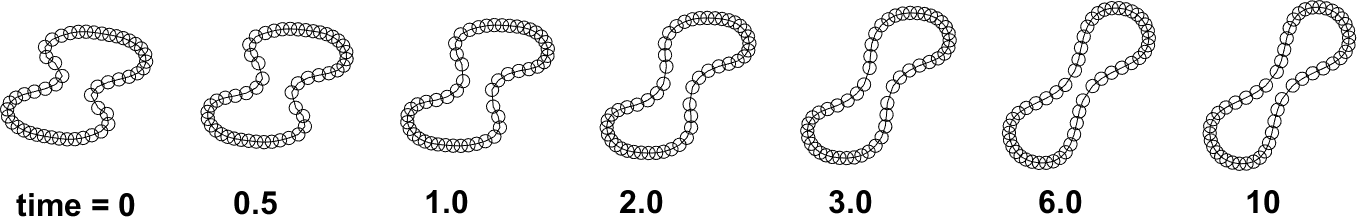}\\
		\hline
	\end{tabular}
	\caption{Cell mid-plane profiles corresponding to Figure \ref{fig:deform}}
	\label{fig:profile}
\end{figure}

Figure \ref{fig:deform} shows deformation snapshots of the evolving cell for the different flow types with spherical harmonic degree $N =24$. The color indicates the norm of the principal membrane stress. Figure \ref{fig:profile} presents the profile snapshots of the cell cross section in the $xz$ plane corresponding to the deformed shapes in Figure \ref{fig:deform}. The behaviors seen in these figures are well documented in the literature \cite{Pozrikidis03, cordasco17}. From Figures \ref{fig:deform}(a, b) and \ref{fig:profile}(a, b), 
it can be seen that the cell rotates clockwise, and is stretched or compressed depending on the angle of inclination. Compared to case (a), the cell in case (b) has larger deformation due to higher shear rate,
and the period of rotation is much higher as well. A star-shape marker serves as tracker of the cell rotating.
The reversed S-shape of the cell mid-plane profile shown in Figure \ref{fig:profile}(b) is also reported by Pozrikidis \cite{Pozrikidis03}. In Figures \ref{fig:deform}(c) and \ref{fig:profile}(c), the cell is stretched at the middle as it translates in the $x$ direction under parabolic flow. Due to the initial inclination, the upper half of the cell moves faster and drags the lower half forward. In the relaxation case, the cell recovers from the highly deformed shape to reach the stress-free state in ambient fluid, as shown in Figure \ref{fig:deform}(d) and \ref{fig:profile}(d).

In the following, various numerical performance aspects of the implicit time integration scheme are discussed:
\begin{enumerate}[label=(\roman*)]
	\item  \textbf{Time step size}: 
	The time-step sizes used by the implicit scheme and the explicit scheme for different spatial resolutions and flow types are reported in Table \ref{tab:timesteps}. Time-step size for implicit method is nearly independent of the spherical harmonic degree, 
	and only depends on the flow type and capillary number. 
	A similar observation was made in \cite{Veerapaneni11} 
	for simulation of a vesicle in simple shear flow with a semi-implicit scheme. 
	\begin{table}
		\centering
		\resizebox{\columnwidth}{!}{
			\begin{tabular}{c|ccc|cccc}
				\hline
				& \multicolumn{3}{c|}{Explicit} &\multicolumn{4}{c}{Implicit}       \\ \hline
				Flow type & (a) & (b) & (c)  
				&(a) & (b)   & (c) & (d) \\ \hline
				$N = 8$   & 1.13E-02 & 1.16E-02 & 1.13E-02 
				& 7.34E-02 & 1.60E-02 & 1.22E-01     & 5.14E-02             \\
				16  & 5.07E-03 & 5.63E-03 & 5.38E-03 
				& 5.10E-02 & 2.26E-02 & 1.17E-01     & 6.55E-02\\
				24  & 3.07E-03 & 3.07E-03 & 3.46E-03 
				& 6.60E-02 & 2.45E-02 & 8.93E-02     & 5.90E-02\\ \hline
			\end{tabular}
		}
		\caption{Steady state time-step sizes used by the explicit and implicit time integration schemes for
			different flow types and spherical harmonic degrees}
		\label{tab:timesteps}
	\end{table}
	
	\begin{figure}
		\centering
		\subfloat
		{
			\resizebox{0.8\textwidth}{!}{\includegraphics{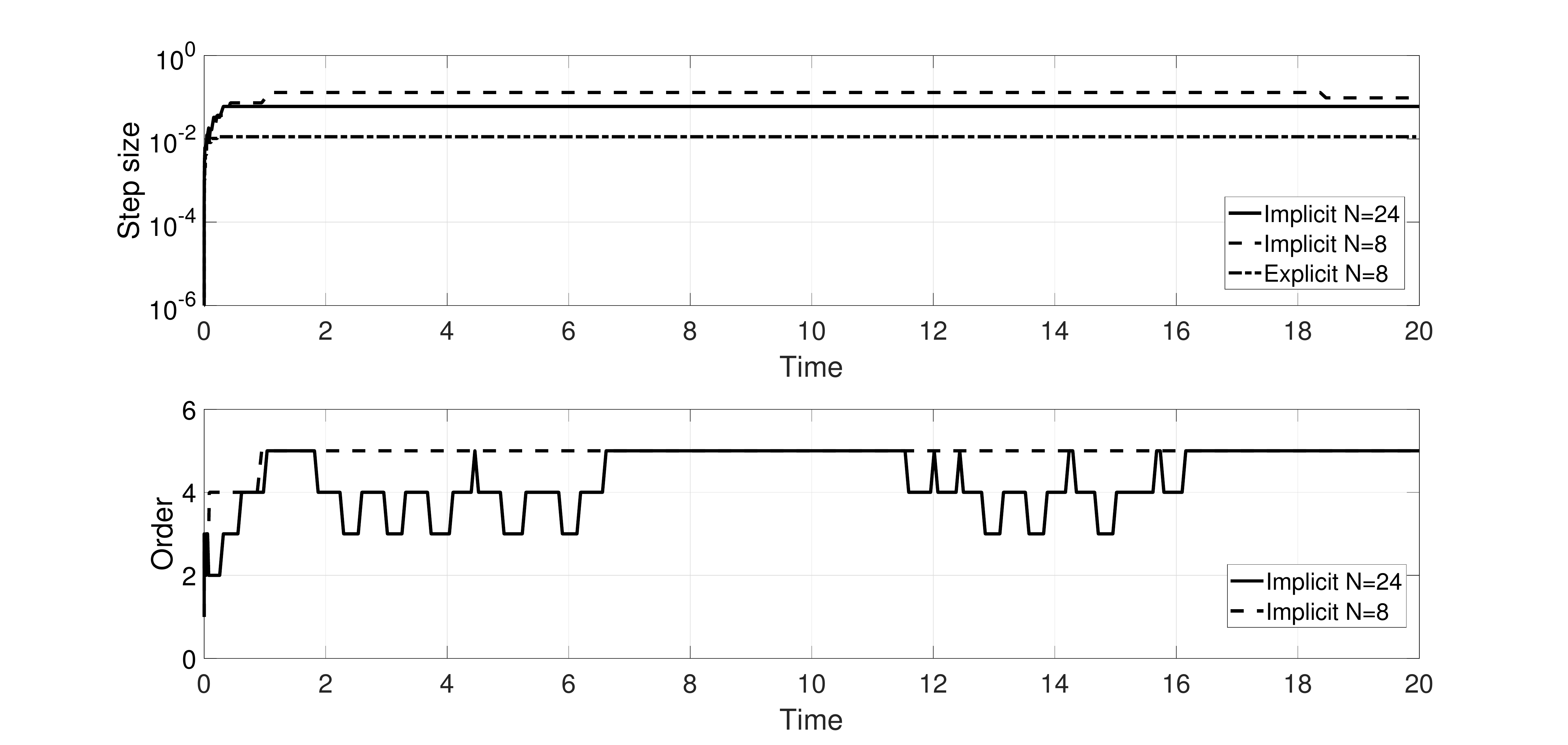}}
		} \hspace{0.05\textwidth}
		\subfloat
		{
			\resizebox{0.8\textwidth}{!}{\includegraphics{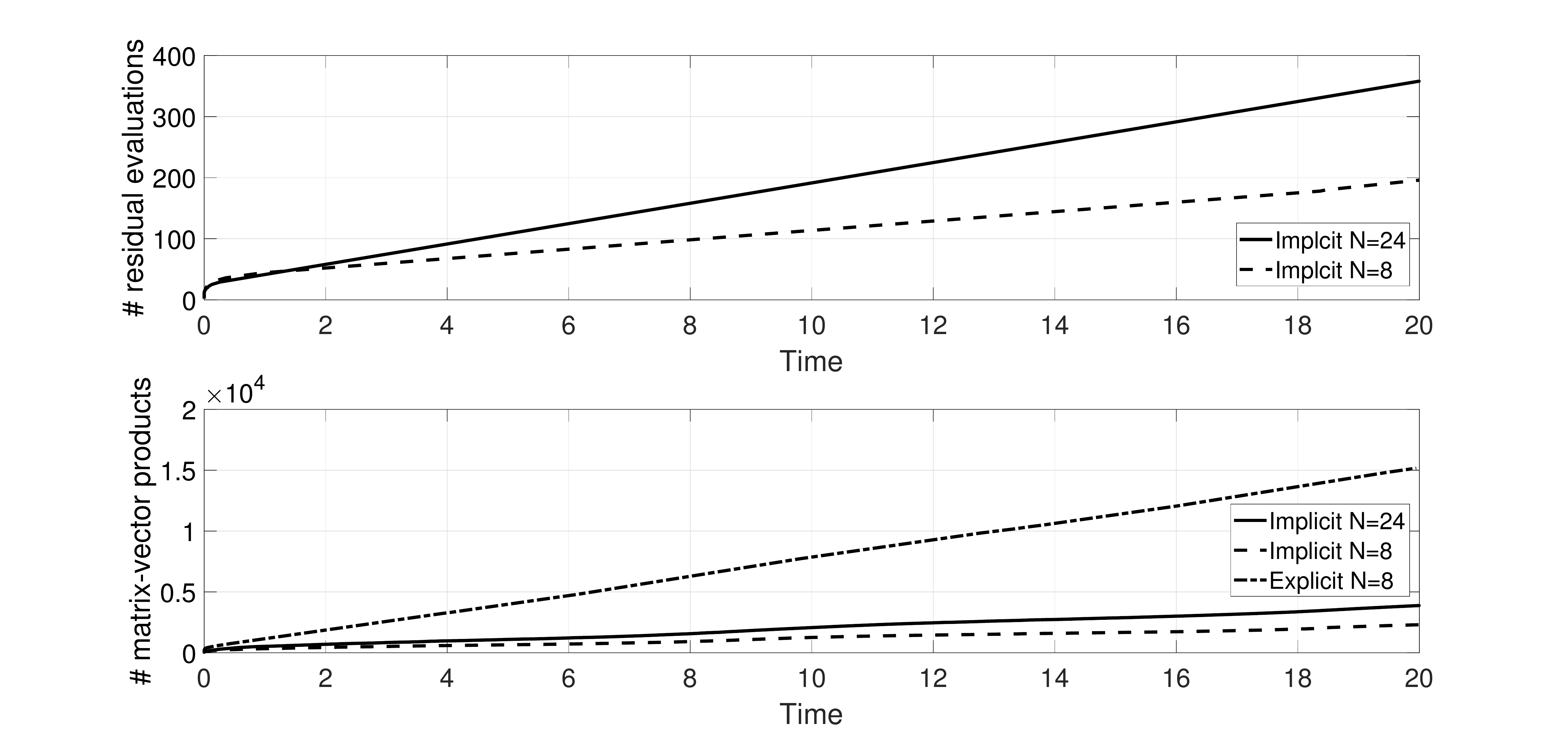}}
		} 
		\caption{\texttt{IDAS} solver statistics for example (a)  
			for $N$ = 24 (solid line), $N$ = 8 (dotted line) and with explicit method for $N$ = 8 (dash-dot line);
			The first subplot shows evolution of time-step size with simulation time --- 
			after adapting the time step, the solver settles at a time step of $6.60 \times10^{-2}$ for $N$ = 24,
			the second subplot shows the the order of backward difference approximation,
			the third and fourth subplots show numbers of residual and Jacobian-vector product evaluations for implicit method and matrix-vector product for explicit method.
			The first and fourth subplots also show corresponding data for explicit integration.}
		\label{fig:stats_Shear}
	\end{figure}
	\begin{figure}
		\centering
		\subfloat
		{
			\resizebox{0.8\textwidth}{!}{\includegraphics{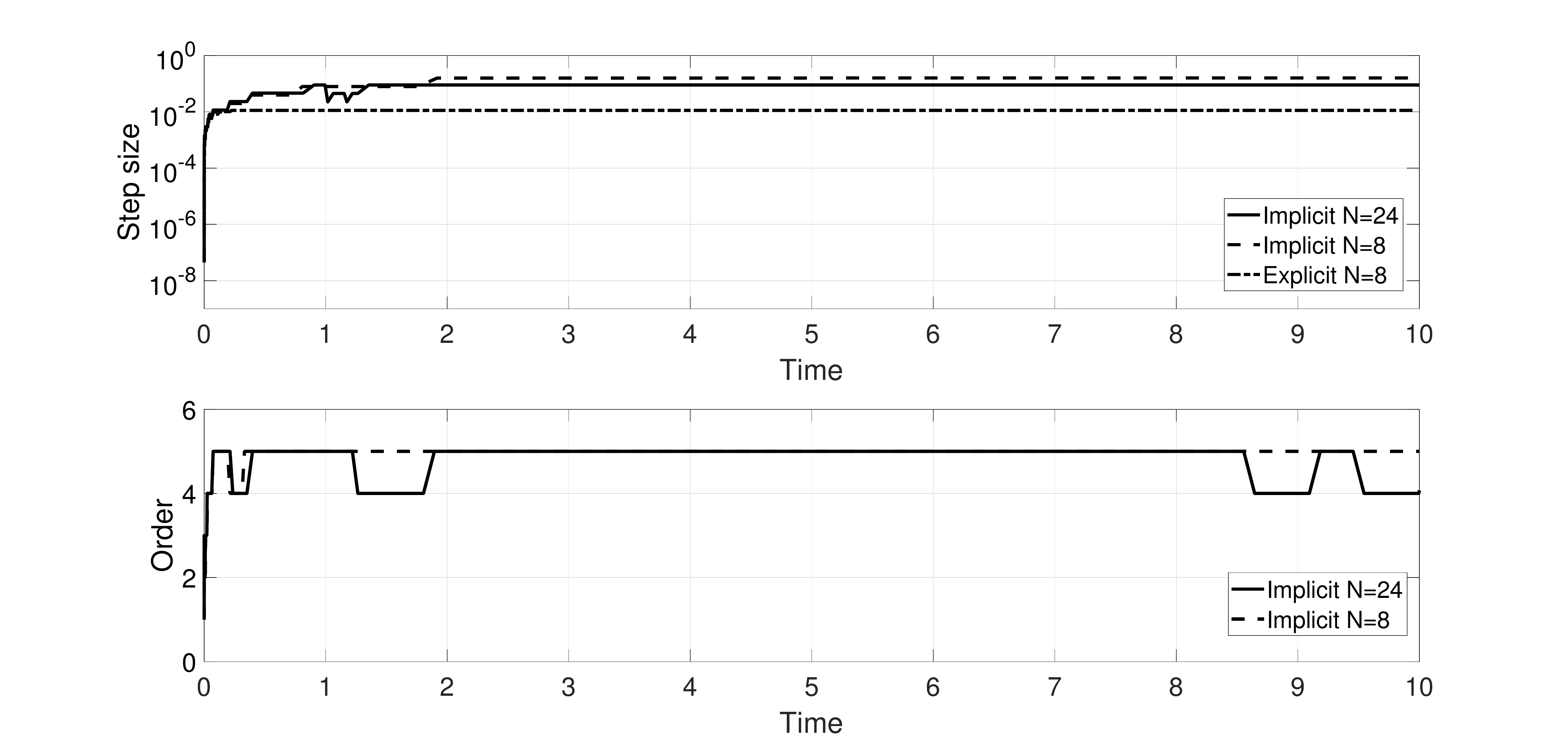}}
		} \hspace{0.05\textwidth}
		\subfloat
		{
			\resizebox{0.8\textwidth}{!}{\includegraphics{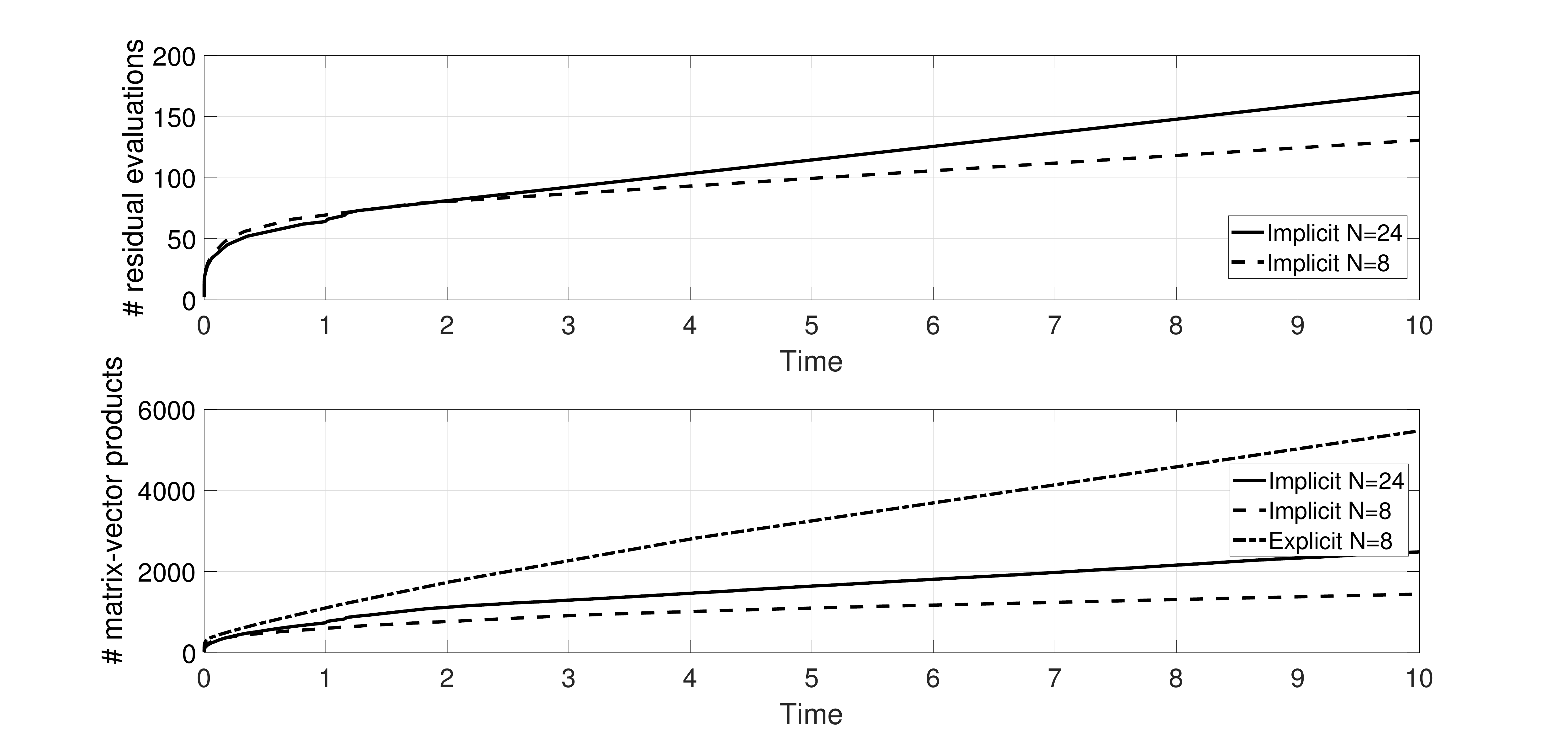}}
		}
		\caption{\texttt{IDAS} solver statistics for example (c) for
			$N$ = 24 (solid line), $N$ = 8 (dotted line) and with explicit method for $N$ = 8 (dash-dot line);
			The first subplot shows evolution of time-step size with simulation time  --- 
			after adapting the time step, the solver settles at a time step of $8.93\times10^{-2}$ for $N$ = 24,
			the second subplot shows the the order of backward difference approximation,
			the third and fourth subplots show numbers of residual and Jacobian-vector product evaluations for implicit method and matrix-vector product for explicit method.
			The first and fourth subplots also show corresponding data for explicit integration.}
		\label{fig:stats_Para}
	\end{figure}
	
	\texttt{IDAS} solver statistics are provided for the shear flow case (a) in Figure \ref{fig:stats_Shear} 
	and for the parabolic flow case (c) in Figure \ref{fig:stats_Para}. 
	In each figure, solid lines show data for $N\!=\!24$, and dotted line for $N\!=\!8$. The data with explicit method for $N\!=\!8$ is also reported as dash-dot line in both first and last subplots for comparison.
	The first subplot shows evolution of time step size during the simulation. 
	The solver adapts the time step based on local truncation error criteria and on convergence of Newton iterations.
	For each example, we set a maximum step size to minimize the number of failed steps.
	The solver settles at a time step of $6.60\!\times\!10^{-2}$ in case (a), and $8.93\!\times\!10^{-2}$ in case (c),
	more than ten times what is possible with an explicit method based on stability considerations alone. 
	
	\item \textbf{Order of backward difference approximation}: 
	The second subplot in each of Figures \ref{fig:stats_Shear} and Figures \ref{fig:stats_Para}
	shows how the order of backward difference approximation employed by the solver. 
	The order is adapted between 1 and 5, again based on local truncation error test and convergence.
	The order of approximation is at least 3, and most often is the highest possible order 5,
	indicating the smoothness of the solution. 
	
	\item \textbf{Number of function and derivative calls needed}: 
	The third and fourth subplots in Figures \ref{fig:stats_Shear} and Figures \ref{fig:stats_Para} 
	report the numbers of residual function evaluations (Figure \ref{fig:IDOE-g})
	and Jacobian matrix-vector multiplications (Figure \ref{fig:IDOE-Dg}).
	Approximately ten Jacobian multiplications are needed for every residual evaluation.
	Since fewer time steps are used, the total number of matrix-vector products
	is comparable to the number of GMRES iterations needed in an explicit method.
	The number of matrix-vector products can potentially be reduced 
	by suitably preconditioning the system. 
	
	\item \textbf{Conservation of different quantities}:
	Besides local truncation error, as other measures of accuracy, 
	we consider conservation of membrane surface area and cell volume, 
	and equilibrium of membrane forces.
	Figure \ref{fig:AV} shows membrane surface area change and volume conservation error
	for the shear case (a). 
	Volume conservation error decreases dramatically with increasing 
	spherical harmonic degree due to the accuracy of representation of the BIE itself.
	The area changes is within $1\%$, indicating the extent to which the penalty parameter, $E_\text{D}$,
	is effective in enforcing membrane inextensibility.
	The largest area changes occur at the times when cell is stretched most by the shear flow 
	when the angle of inclination of the cell is  $\pm 45^\circ$.
	Figure \ref{fig:FM} demonstrates equilibrium of forces and moments on 
	the entire cell.
	
	\item \textbf{Comparison with the explicit integration}
	Due to the variabilities in how time step size can be adapted, 
	tolerances for nonlinear and linear solvers, truncation error control etc., 
	it is not possible to make an exact quantitative comparison between the performances
	of implicit and explicit time integration schemes. 
	Qualitatively, however, it is clear that steady state time-step size used by the implicit time integration scheme depends on the imposed flow type and is not affected by the spherical harmonics order (see Table \ref{tab:timesteps}). On the other hand, the steady time-step size used by the explicit method decreases with increasing spherical harmonic degree and, interestingly, is not affected by the flow strength and type. 
	The time step dictated by stability in an explicit method is approximately 1/10th of the time step size used by 
	\texttt{IDAS} for truncation error control. 
	Comparison of computational effort between the IDAS and an explicit method is 
	made using the number of matrix-vector multiplications (although a Jacobian matrix-vector multiplication is more expensive than the residual matrix-vector multiplication). 
	The number of matrix-vector multiplication performed by IDAS is approximately 1/7th of the number of GMRES iterations in an explicit method as seen in the fourth subplots of Figures \ref{fig:stats_Shear} and \ref{fig:stats_Para}.

	\begin{figure}
		\centering
		\includegraphics[width=\textwidth]{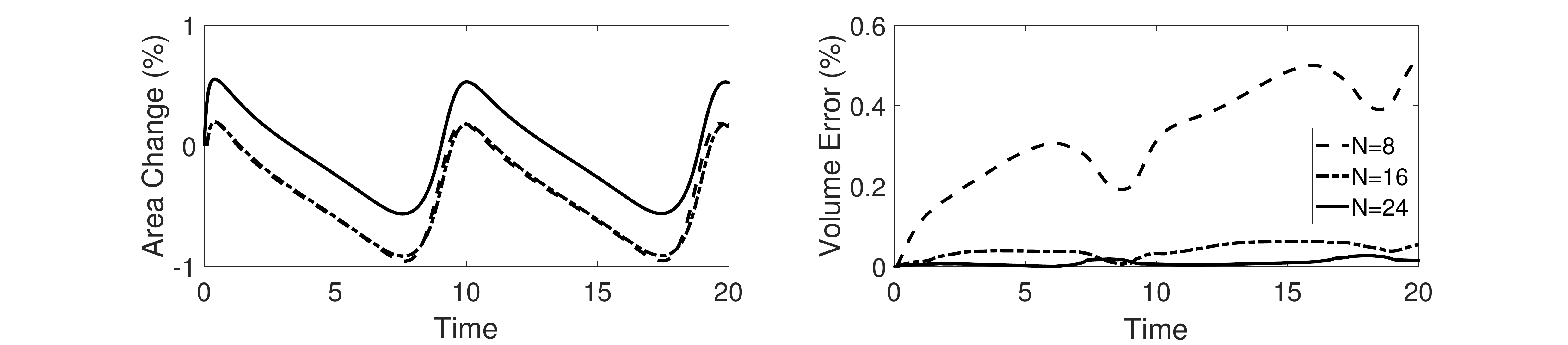}
		\caption{Membrane surface area change and volume conservation error 
			for shear flow example (a) with $k =  100$ $\SI{}{\per\second}$ for different 
			spherical harmonic degrees $N$.}
		\label{fig:AV}
	\end{figure}
	\begin{figure}
		\centering
		\includegraphics[width=0.9\textwidth]{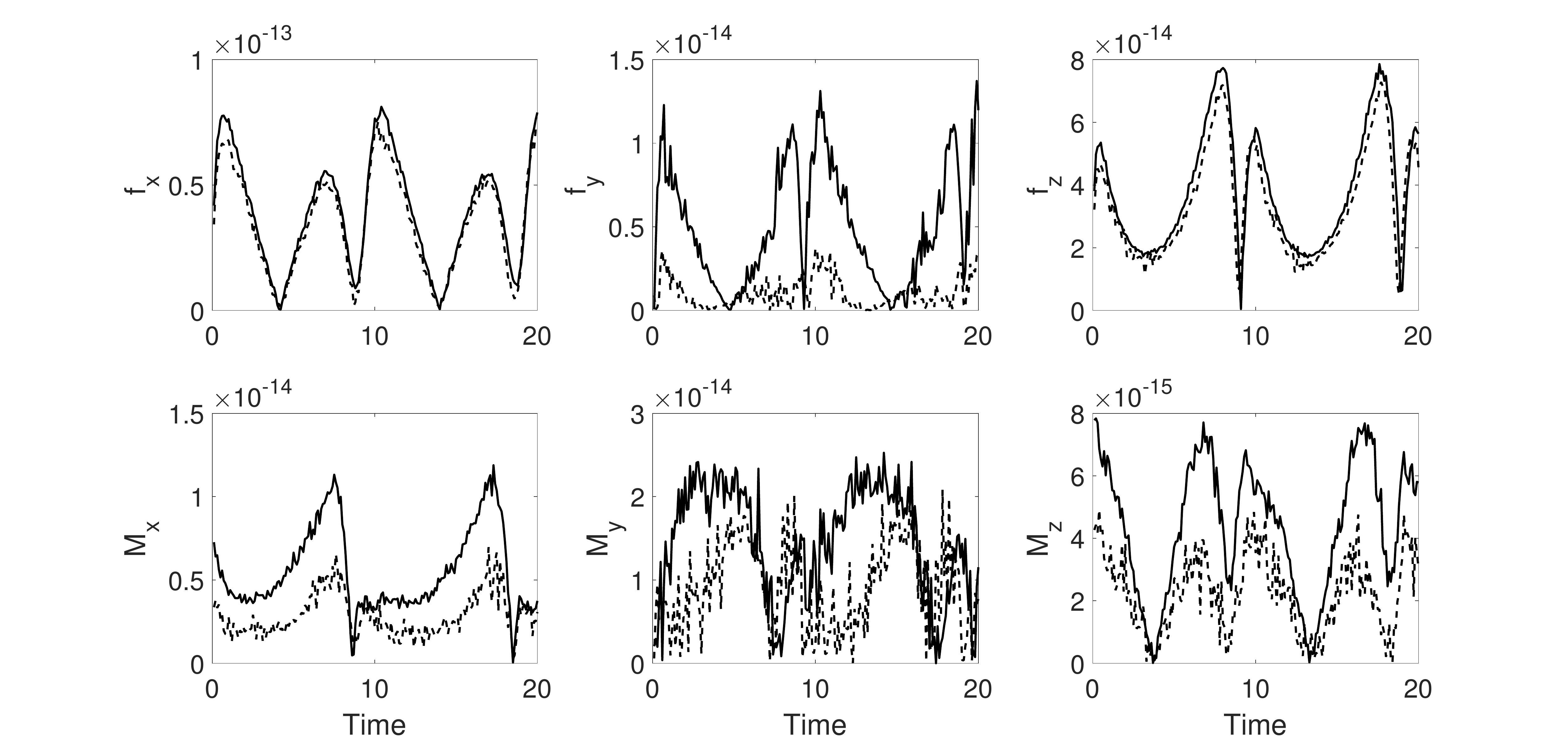}
		\caption{Summation of membrane forces and moments for shear flow example (a) with $k =  100$ $\SI{}{\per\second}$ for $N$ = 24 (solid line) and $N$ = 8 (dotted line).}
		\label{fig:FM}
	\end{figure}
\end{enumerate}

\section{Concluding remarks}\label{sec:Conclusion}
An implicit time integration approach has been presented in connection with
a spectral boundary integral method for computation of red blood cell flow.
The simple problem of cell flow in an ambient fluid (rather than in a vessel)
has been considered to describe of the method.
The flow is represented by a boundary integral equation (BIE) 
expressed in weakly singular form.
Cell configurations are represented by smooth deformations of the unit sphere.
The cell configuration and velocity field are discretized using
spherical harmonic basis functions.
The cell membrane is modeled using elastic thin-shell mechanics,
formulas from which are written in a manner that requires only
the operations of gradient of scalar functions and 
divergence of vector fields on the unit sphere in spherical coordinates.
This BIE representing the flow has the form of an implicit ordinary differential equation (IODE)
in the cell configuration, and solutions are computed using
an off-the-shelf IODE solver based on backward differences formulas.
The solver adapts the time step and order of backward difference approximation
to control local truncation error and to regulate convergence of Newton's method.
Newton's method employed in implicit time integration necessitates
computation of Jacobians of the IODE with respect to the configuration and the velocity field.
Rather than assemble such Jacobians, Jacobian-vector products, i.e. directional derivatives
are calculated. This calculation, which itself is the evaluation of a different BIE,
is aided by the expression of the flow BIE in weakly singular format.
Numerical examples have been presented to demonstrate the working of the implicit approach.
They illustrate that larger time steps can be used compared to explicit methods,
where time-step size is limited by stability considerations.
The number of Jacobian-vector products is comparable to the number of 
GMRES iterations of the flow BIE in explicit methods.
Source code for the implementation of the implicit method 
is included online supplemental material that accompany this note.

Much further development is needed to employ this implicit method in
broader problems of interest.
Suitable preconditions must be constructed for the directional derivative BIEs
to reduce the number of Jacobian-vector products further, and hence computational effort.
Fast summation methods must also be implemented for the derivative BIEs.
Most importantly, the method must be extended for the presence of 
other geometries such as vessels and slits.
While these do not alter the mathematical formulation significantly,
they raise the possibility of interpenetration when large time steps are used,
In such cases, interpenetration must be flagged to the solver.
Closest-point calculations used in evaluation of near-singular 
integrals can potentially be adapted to detect interpenetration.
These are topics of current inquiry.

\bibliographystyle{plain}
\bibliography{FullyImplicitSpectralBoundaryIntegralComputationOfRBCFlow}  
\end{document}